\documentclass[11pt,reqno]{amsart}

\usepackage{euscript}
\usepackage{pb-diagram}
\usepackage{lamsarrow}
\usepackage{pb-lams}
\usepackage{epsfig}

\theoremstyle{plain}

\theoremstyle{definition}

\DeclareMathOperator{\Hom}{{\rm Hom}}

\newcommand\bC{{\mathbb C}}

\newcommand{\bQ}{{{\mathbb Q}}}
\newcommand\bR{{\mathbb R}}

\newcommand\bZ{{\mathbb Z}}

\newcommand{\Z}{{\mathfrak Z}}

\newcommand\aug{\mathfrak{aug}}

\newcommand{\p}{{\mathcal P}}

\newcommand{\R}{\mathcal R}

\newcommand{\C}{{\mathcal C}}

\newcommand{\ra}{\rightarrow}

\def\inpr{\mathbin{\hbox to 6pt{\vrule height0.4pt width5pt depth0pt \kern-.4pt \vrule height6pt width0.4pt depth0pt\hss}}}

\newcommand{\si}{{\sigma}}

\newcommand{\eps}{\epsilon}
\newcommand{\ssw}{{\bf sw}}

\newcommand{\et}{\EuScript{T}}

\newcommand{\bms}{\mbox{\boldmath$s$}}

\newcommand{\x}{o}

\def\mmod{\mbox{mod}}
\let\d\partial

\def\C{\mathbb C}
\def\Q{\mathbb Q}
\def\R{\mathbb R}
\def\Z{\mathbb Z}

\begin{document}

\title{Seiberg-Witten invariants and surface singularities II\\
(singularities with good $\C^*$-action)}

\author{Andr\'as N\'emethi}
\address{Department of Mathematics\\Ohio State University\\Columbus, OH 43210}
\email{nemethi@math.ohio-state.edu}
\urladdr{http://www.math.ohio-state.edu/ \textasciitilde nemethi/}
\author{Liviu I. Nicolaescu}
\address{University of Notre Dame\\Notre Dame, IN 46556}
\email{nicolaescu.1@nd.edu} \urladdr{http://www.nd.edu/
\textasciitilde lnicolae/}
\thanks{The first author is partially supported by NSF grant DMS-0088950;
the second author is partially supported by NSF grant
DMS-0071820.}



\begin{abstract} We verify the conjecture formulated in \cite{NN} for any
normal surface singularity which admits a good $\C^*$-action. The main result
connects the Seiberg-Witten invariant of the link
(associated with a certain ``canonical'' $spin^c$ structure) with the
geometric genus of the singularity.

As a by-product, we compute the Seiberg-Witten monopoles of the link
(associated with the canonical $spin^c$ structure  and the natural Thurston
 metric) in terms of its Seifert invariants.
Additionally, we  also determine  in terms of the Seifert invariants
  all the Reidemeister-Turaev sign-refined
torsion of the link (associated with  any $spin^c$ structure).
\end{abstract}

\maketitle
\pagestyle{myheadings}
\markboth{{\normalsize Andr\'as N\'emethi and Liviu I. Nicolaescu}}{
{\normalsize Seiberg-Witten invariants and surface singularities II}}

\section{Introduction}

The present article is a natural continuation of \cite{NN},
where the authors formulated  a very
general conjecture  which relates the topological and the
analytical invariants  of a complex normal surface singularity
whose link is  a rational homology sphere.

Let $(X,0)$ be a normal two-dimensional analytic  singularity. It
is well-known that from a topological point of view, it   is
completely characterized by its link $M$, which is an oriented
3-manifold. Moreover, by a result of Neumann \cite{NP}, any
decorated resolution graph of $(X,0)$ carries the same information
as $M$. A property of $(X,0)$ will be called \emph{topological} if
it can be determined from  $M$, or equivalently, from any
resolution graph of $(X,0)$.
For example, for a given resolution, if we take
the canonical divisor $K$, and the number $\# {\mathcal V}$
of irreducible components of the exceptional divisor of the resolution,
then $K^2+\#{\mathcal V}$ is independent of the choice of the
resolution, it is an invariant of the link $M$ (cf. \ref{ss: 2.5}).

Our interest is to investigate the possibility to express the
geometric genus (which is, by its very definition, an analytic invariant of the
singularity) in terms of topological invariants of the link.

Let us start with a brief historical survey.
M. Artin  proved in \cite{Artin,Artin2}
that the rational singularities (i.e. $p_g=0$) can be
characterized  completely from the graph.
In \cite{Lauferme}, H. Laufer   extended  Artin's results to minimally elliptic
singularities, showing that Gorenstein singularities with $p_g=1$ can be
characterized topologically.
Additionally, he noticed that the program breaks for
more complicated singularities (see also the comments in \cite{NN} and
\cite{Neminv}).
On the other hand, the first author noticed in \cite{Neminv} that
Laufer's counterexamples do not signal  the end of the program. He
conjectured  that if we restrict ourselves to the case of those
Gorenstein singularities whose links are rational homology spheres
then  $p_g$ is topological. This was carried out explicitly for elliptic
singularities in \cite{Neminv} (partially based on some results of
S. S.-T. Yau, cf.  e.g.  with \cite{Yau}).

For Gorenstein singularities, in the presence of a smoothing
with Milnor fiber $F$, the above question can be  reformulated in terms of
the signature $\si(F)$ and/or  the
topological Euler characteristic $\chi_{top}(F)$ of $F$. Indeed,
via some results of Laufer, Durfee, Wahl  and
Steenbrink, for Gorenstein singularities, any of $p_g$,
$\si(F)$ and $\chi_{top}(F)$ determines  the remaining two  modulo
$K^2+\#{\mathcal V}$ (see e.g.  \cite{LW}).

This fact creates the bridge connecting
 the above problem with the following list of results
about the signature $\si(F)$.
Fintushel and Stern proved in \cite{FS}  that
for a hypersurface Brieskorn  singularity whose link is an {\em
integral} homology sphere, the Casson invariant $\lambda(M)$  of
the link $M$ equals $-\si(F)/8$.  This fact was generalized
by Neumann  and Wahl in \cite{NW}. They proved the same statement
for all Brieskorn-Hamm complete intersections and suspensions of plane
curve singularities (with the same assumption about the link).
Moreover, they conjectured the validity
of the formula for any isolated complete intersection  singularity
(with the same restriction about the link). For some other
conjectures about these singularities, the reader can also consult
\cite{NW2}.

The goal of \cite{NN} was to generalize this conjecture for
smoothing of Gorenstein singularities with rational homology sphere link.
In fact, more generally (i.e. even if the singularity is not smoothable)
\cite{NN}
conjecturally  introduced an  {\em ``optimal'' topological upper bound}
for $p_g$ in the following sense:
it is a topological upper bound for $p_g$ for {\em any} normal
surface singularity, but additionally,
for Gorenstein singularities it yields exactly $p_g$.
(Such an ``optimal'' topological upper bound
for {\em elliptic} singularities  is
the length of the elliptic sequence, introduced and studied by  S. S.-T. Yau,
see e.g. \cite{Yau},  and Laufer.)

The conjecture in \cite{NN}  replaces  the
Casson invariant $\lambda(M)$ by a certain Seiberg-Witten
invariant of the link, i.e. by the sum of the Casson-Walker
invariant and a certain Reidemeister-Turaev sign-refined torsion
invariant.

We recall (for details, see the first part \cite{NN}, and the references
listed there) that the set of
Seiberg-Witten invariants associates with any $spin^c$ structure
$\si$ of $M$ a rational number $\ssw^0_M(\si)$. In \cite{NN}
 we introduced a  ``canonical'' $spin^c$
structure $\si_{can}$  of  $M$. This can be done as follows. The
(almost) complex structure on $X\setminus\{0\}$  induces a natural
$spin^c$ structure on $X\setminus \{0\}$. Then $\si_{can}$, by definition,
is its restriction to $M$.
(An equivalent definition can be done using a canonical quadratic
function constructed by Looijenga and Wahl in \cite{LW}.)
The point is that this structure
 depends only on the topology  of  $M$ alone.

\vspace{2mm}

\noindent We  are now ready to recall the  conjecture from \cite{NN}.

\subsection{Conjecture.}\label{MC}\cite{NN} {\em Assume that $(X,0)$ is a
normal surface singularity whose link $M$ is a rational homology
sphere. Let $\si_{can}$ be the canonical $spin^c$  structure on
$M$. Then, conjecturally, the following facts hold.

\vspace{2mm}

\noindent (1) For any $(X,0)$, there is a topological upper bound
for $p_g$ given by:
$$\ssw^0_M(\si_{can})-\frac{K^2+\#{\mathcal V}}{8}\geq p_g.$$
(2) If $(X,0)$ is rational or Gorenstein, then in (1) one has equality.\\
(3) In particular, if $({\mathcal X},0)$ is a smoothing of a
Gorenstein singularity $(X,0)$ with Milnor fiber $F$, then
$$\ssw^0_M(\si_{can})=-\frac{\si(F)}{8}.$$}
[Notice that if $(X,0)$ is numerically Gorenstein and $M$ is a $\Z_2$--homology
sphere then $\si_{can}$ is the unique spin structure of
$M$; if $M$ is an integral  homology sphere then in the above
formulae $\ssw^0_M(\si_{can})=\lambda(M)$, the Casson invariant
of $M$.]

\vspace{3mm}

In \cite{NN}, the conjecture was verified for cyclic quotient and
du Val singularities, Brieskorn-Hamm complete
intersections, and some rational and minimally elliptic singularities.

In \cite{NN} Remark 4.6(2),  we  noticed that it is very likely that parts
(2) and (3)  of the conjecture
are  true even for a larger class of singularities, e.g.
for $\Q$-Gorenstein singularities, or some families of singularities with
some kind of additional analytical rigidity.

In fact, this is exactly the case for normal surface singularities with
a (good) $\C^*$-action. It is well-known that a complex
 affine algebraic variety  $X$ admits a $\C^*$-action
if and only if the affine coordinate ring $A$ admits a grading
$A=\oplus_kA_k$. Following Orlik-Wagreich, we say that the action is
{\em good} if $A_k=0$ for $k<0$ and $A_0=\C$. This means that the point
$0$ corresponding to the maximal ideal $\oplus_{k>0}A_k$ is the only
fixed point of the action. Additionally, we will assume that
$(X,0)$ is normal.

For these singularities, in this article we prove:

\vspace{2mm}

\subsection{Theorem.}\label{Th1} \
 {\em Let $(X,0)$ be a normal surface singularity
with a good  $\C^*$-action. Then }
$$\ssw^0_M(\si_{can})-\frac{K^2+\#{\mathcal V}}{8}=p_g.$$

\vspace{2mm}

The proof is based  in part on Pinkham's formula \cite{Pinkh} of $p_g$
in terms of the Seifert invariants of the link (cf. also
with Dolgachev's work about weighted homogeneous singularities;
see  e.g. \cite{Dolg}).
On the other hand,  in the proof we use the formulae for
$K^2+\#{\mathcal V}$ and the Reidemeister-Turaev torsion determined
in \cite{NN}, and  a formula for the Casson-Walker invariant proved in
\cite{Lescop}. The Reidemeister-Turaev torsion formally shows big
similarities with the Poincar\'e series associated with
the graded coordinate ring of the universal  abelian cover of $(X,0)$.
In the proof we also borrowed some technique of Neumann  applied by him
for this Poincar\'e series \cite{Neu}.

Notice that in the above theorem we do not require for $(X,0)$ to be
Gorenstein. (This  is replaced by the existence of the $\C^*$-action.)
On the other hand, the theorem has the following corollary
(which can also be  applied  for singularities without $\C^*$-action).

\vspace{2mm}

\subsection{Corollary.} \label{C1} {\em Assume that the link  of a
normal surface singularity $(X,0)$
is a rational homology sphere Seifert 3-manifold.
If $(X,0)$ is rational, or minimally elliptic, or Gorenstein elliptic,
then the statements of the above conjecture are true for  $(X,0)$.
(I.e., $(X,0)$ satisfies (1) with equality; and also (3),  provided that
the  additional assumptions of (3) are satisfied.)}

\vspace{2mm}

\noindent Indeed, in the case of these singularities, all the numerical
invariants involved in the conjecture are characterized by the link.
Moreover, each family contains a special representative which admits a
good $\C^*$-action. In fact, the above corollary can
automatically  be extended to  any
family of singularities with these two properties.

\vspace{2mm}

The paper is organized as follows. In section 2 we review the needed
definitions and results.  For a more complete picture and list of references
the reader is invited to
consult \cite{NN}. Section 3 starts with a  theorem
(cf. \ref{m1}) which connects four topological numerical invariants
of the link. These invariants are: the Reidemeister-Turaev torsion, the
Casson-Walker invariant, the Dolgachev-Pinkham invariant $DP_M$
(which is the topological candidate for $p_g$), and finally
$K^2+\#{\mathcal V}$ (which can be identified with the Gompf invariant,
cf. \ref{mlb}).
This result implies the above theorem via Pinkham's result \cite{Pinkh}
(cf. \ref{ss: inv}(10)).

Finally, we have included at  the end of section 3
a complete and explicit  description
the Reidemeister-Turaev torsion $\et_{M,\si}(1) $ for
{\em any} $spin^c$ structure $\si$ in terms of the Seifert invariants of $M$
(the computation follows Neumann's method mentioned above).

\subsection{Remark.}\label{C2} Theorem \ref{m1} has a
consequence  interesting from the point of gauge theory as well. Recall that
the modified Seiberg-Witten invariant is the sum
$\ssw^0_M(\si)=\ssw_M(\si,u)+KS_M(\si,u)/8$ of the Seiberg-Witten
monopoles and ($1/8$ of)  the Kreck-Stolz invariant (cf. \ref{ss:
3.2}). The $(\si_{can},u)$-monopoles were described in
\cite{MOY,Nico2}, and the Kreck-Stolz invariant was described  in
\cite{Nico2}. Nevertheless, their explicit computation in terms of
the Seifert invariants in the general case meets some
difficulties. Our theorem \ref{m1} implies  the following.

\smallskip

Assume that $M$ is a rational
homology sphere Seifert 3-manifold with $e<0$.
Then, for any \emph{good} parameter $u$,
 the signed number of Seiberg-Witten $u$-monopoles can be computed as follows:
$$\ssw_M(\si_{can},u)=-\frac{KS_M(\si_{can},u)}{8}+\frac{K^2+\#
{\mathcal V}}{8}+DP_M.$$
Recall that the Seifert manifold $M$ admits a natural metric $g_0$,
the so called Thurston metric. If  $u=(g_0,0)$,
the invariants appearing on  the right hand side have very explicit
expression in terms of the Seifert  invariants of $M$. Indeed,
$KS_M(\si_{can},g_0,0)$ is determined in \cite{Nico2} (see also \cite{NN}),
$K^2+\#{\mathcal V}$ is determined in \cite{NN} (see
\ref{ss: inv}(6) here), and $DP_M$ was introduced in \cite{Pinkh}
(it is \ref{ss: inv}(9) here); in fact the last two are $u$-independent.
In particular, when the parameter $(g_0,0)$ is good, the above identity
determines $\ssw_M(\si_{can},g_0,0)$
explicitly  in terms of the Seifert invariants.

For results in this direction, see \cite{MOY,Nico2}.

\medskip

\vspace{2mm}

\section{Preliminaries}

\subsection{Definitions.}\label{ss: 4.1} Let $(X,0)$ be a normal surface
singularity. Consider the  holomorphic line bundle
$\Omega^2_{X\setminus\{0\}}$ of holomorphic 2-forms on $X\setminus
\{0\}$. If this line bundle is holomorphically trivial then we say
that $(X,0)$ is \emph{Gorenstein}. Let $\pi:\tilde{X}\to X$
be a resolution over a sufficiently small Stein representative  $X$ of the germ
$(X,0)$. Then $p_g(X,0):=\dim \, H^1(\tilde{X}, {\mathcal
O}_{\tilde{X}})$ is finite and independent of the choice of $\pi$.
It is called the {\em geometric genus } of $(X,0)$.

\subsection{The link and its canonical $spin^c$ structure.} \label{ss: 2.1}
Let $(X,0)$ be a normal surface singularity
embedded in  $(\C^N,0)$. Then for $\epsilon$ sufficiently small
the intersection $M:=X\cap S_{\epsilon}^{2N-1}$ of a
representative $X$ of the germ with the sphere
$S_{\epsilon}^{2N-1}$ (of radius $\eps$) is a compact oriented
3-manifold, whose oriented $C^{\infty}$ type does not depend on
the choice of the embedding and $\epsilon$. It is called the link
of $(X,0)$.
In this article we will assume that $M$ is a rational homology
sphere, and we write $H:=H_1(M,\Z)$.

The almost  complex structure on $X\setminus\{0\}$
determines a $spin^c$ structure on $X\setminus \{0\}$,
whose restriction to $M$ will be denoted by  $\si_{can}\in Spin^c(M)$.
It turns out that $\si_{can}$ depends only on the oriented $C^{\infty}$
type of $M$ (for details  see \cite{NN}, cf. also with \cite{LW}).

\subsection{The Seiberg-Witten invariants of $M$.}\label{ss: 3.2}
To describe the
Seiberg-Witten invariants one has to consider  an additional geometric
data belonging to  the space of  parameters
\[
\p=\{ u=(g,\eta);\quad g=\mbox{Riemann metric},\;\;\eta=\mbox{closed
two-form} \}.
\]
Then for each $spin^c$ structure $\si$ on $M$  one defines the
$(\si,g,\eta)$-{\em Seiberg-Witten monopoles}.
For a generic parameter $u$, the Seiberg-Witten invariant  $\ssw_M(\si,u)$ is  the signed monopole count.
This integer depends on the choice of the parameter $u$ and thus
it is not a topological  invariant.   To obtain an invariant of $M$,  one
needs  to  alter this monopole count. The needed additional
contribution is the {\em Kreck-Stolz} invariant $KS_M(\si,u)$
(associated with the data $(\si,u)$),
cf. \cite{Lim} (or see \cite{KrSt} for the original ``spin version'').
Then, by \cite{Chen1,Lim,MW}, the rational number
\[
\frac{1}{8}KS_M(\si, u)+\ssw_M(\si, u)
\]
is independent   of $u$ and thus  it is a topological invariant of
the pair $(M,\si)$. We denote this {\em modified Seiberg-Witten
invariant}  by \ $\ssw_M^0(\si)$.

\subsection{The Reidemeister-Turaev torsion and the
Casson-Walker invariant.} \label{ss: 3.4}
For any $spin^c$ structure $\si$ on $M$, we denote by
$$\et_{M,\si}=\sum_{h\in H}\et_{M,\si}(h)\, h\in {\bQ}[H]$$
 the sign refined
{\em Reidemeister-Turaev torsion} associated with $\si$
(for its detailed description, see \cite{Tu5}).
It is convenient to think of $\et_{M,\si}$ as a function
$H\ra {\bQ}$ given by $h\mapsto \et_{M,\si}(h)$.
The  {\em augmentation map} $\aug: {\bQ}[H]\ra {\bQ}$
is defined by $\sum a_h\, h\mapsto \sum a_h$. It is known that
$\aug(\et_{M, \si})=0$.

Denote by $\lambda(M)$ the Casson-Walker invariant of $M$
normalized as in \cite[\S 4.7]{Lescop}. Then by  a result of the second author
\cite{Nico5}, one has:
\begin{equation*}
\ssw^0_M(\si)=\frac{1}{|H|}\lambda(M)+\et_{M,\si}(1).
\tag{1}
\label{3.7.2}
\end{equation*}
Below we will present a formula for
$\et_{M,\si}$ in terms of Fourier transform. For this, consider
the Pontryagin dual $\hat{H}:=\Hom(H, U(1))$  of $H$.
Then a function $f: H \ra {\bC}$ and
its Fourier transform $\hat{f}:\hat{H}\ra {\bC}$ satisfy:
 \[
\quad \hat{f}(\chi)=\sum_{h\in H} f(h)\bar{\chi}(h);\ \
 f(h)=\frac{1}{|H|}\sum_{\chi\in \hat{H}}\hat{f}(\chi)\chi(h).
 \]
 Notice  that $\hat{f}(1)=\aug(f)$,
in particular $\hat{\et}_{M,\si}(1)=\aug(\et_{M, \si})=0$.

\subsection{$M$ as a plumbing manifold.} \label{ss: 2.5}
Fix a sufficiently small (Stein) representative  $X$ of $(X,0)$
and let $\pi:\tilde{X}\to X$
be a resolution of the singular point $0\in X$. In particular,
$\tilde{X}$ is smooth, and $\pi$ is a biholomorphic  isomorphism
above $X\setminus \{0\}$. We will assume that the exceptional
divisor $E:=\pi^{-1}(0)$ is a normal crossing divisor with
irreducible components $\{E_v\}_{v\in {\mathcal V}}$. Let
$\Gamma(\pi)$ be the  dual resolution graph associated with $\pi$
decorated with the self intersection numbers $\{E_v\cdot E_v\}_v$.
$\Gamma(\pi)$ can be identified with a plumbing graph, and $M$
with a plumbing 3-manifold constructed from $\Gamma(\pi)$ as its
plumbing graph.
Since $M$ is a rational homology sphere, all the irreducible
components $E_v$ of $E$ are rational, and $\Gamma(\pi)$ is a tree.

Let $D_v$ be a  small transversal disc to $E_v$.
In fact $\d D_v$ can be considered as the generic
fiber of the $S^1$-bundle over $E_v$ used in the plumbing
construction of $M$. Consider the elements
$g_v:=[\d D_v]$ ($v\in {\mathcal V}$) in $H$. It is not difficult to see
that they, in fact,  generate $H$.

For the degree of any vertex $v$  (i.e. for $\#\{w: E_w\cdot E_v=1\}$) we
will use the notation $\delta_v$.

Next, we define the \emph{canonical cycle} $Z_K$ of $(X,0)$ associated with
the resolution $\pi$. This is a {\em rational}  cycle
$Z_K=\sum_{v\in {\mathcal V}}r_vE_v$, $r_v\in\Q$,
 supported  by the exceptional
divisor $E$, and  defined by (the adjunction formula):
\begin{equation*}
Z_K\cdot E_v=E_v\cdot E_v+2 \ \ \ \mbox{for any $v\in {\mathcal V}$}.
\label{2.5*}
\end{equation*}
Since the intersection matrix $\{E_v\cdot E_w\}_{v,w}$ is nondegenerate,
the above equation  has a unique solution.
$(X,0)$ is called {\em numerically Gorenstein}
if $r_v\in \Z$ for each $v\in {\mathcal V}$.

The rational number $Z_K\cdot Z_K$ will be denoted by $K^2$.
Let $\#{\mathcal V}$ denote the number of irreducible components
of $E=\pi^{-1}(0)$. Then $K^2+\#{\mathcal V}$ does not depend on
the choice  of the resolution $\pi$, it is an invariant of $M$.

\vspace{2mm}

\noindent
The main object of this paper is a normal surface singularity
$(X,0)$ with a {good} $\C^*$--action. It is well--known that the link of such
a singularity is a Seifert 3-manifold, and the minimal resolution graph
is a star-shaped graphs. In
these case it is convenient to express all the topological
invariants of $M$ in terms of their
Seifert invariants. In the next subsections  we list
briefly the definitions, notations and some of the needed properties.

\subsection{The Seifert invariants.}\label{ss: s1} \cite{JN,Neu,NeR}\
Consider a Seifert fibration $\pi:M\to \Sigma$. In
our situation, since $M$ is a rational homology sphere, the base
space $\Sigma$ has genus zero.

Consider a set of points $\{x_i\}_{i=1}^{\nu}$ in such a way that
the set of fibers $\{\pi^{-1}(x_i)\}_i$ contains the set of
singular fibers. Set $O_i:= \pi^{-1}(x_i)$. Let $D_i$ be a
small disc  in $X$ containing  $x_i$,
$\Sigma':=\Sigma\setminus \cup_iD_i$ and $M':=\pi^{-1}(\Sigma')$.  Now,
$\pi:M'\to \Sigma'$ admits sections, let $s:\Sigma'\to M'$ be one of them.
Let $Q_i:=s(\d D_i)$ and let $H_i$  be  a circle fiber in
$\pi^{-1}(\d D_i)$. Then in $H_1(\pi^{-1}(D_i),\Z)$ one has
$H_i\sim \alpha_iO_i$ and $Q_i\sim -\beta_iO_i$ for some integers
$\alpha_i>0 $ and $
 \beta_i$ with $(\alpha_i,\beta_i)=1
$. The set $((\alpha_i,\beta_i)_{i=1}^\nu)$ constitute the set of
\emph{(unnormalized) Seifert invariants}. The number
\[
e:=- \sum\frac{\beta_i}{\alpha_i}
\]
 is called the \emph{(orbifold) Euler number} of $M$.
$M$ is a link of singularity if and only if  $e<0$.

Replacing the section by another one, a different choice  changes
each $\beta_i$ within its residue class modulo $\alpha_i$ in such
a way that the sum $e=-\sum_i(\beta_i/\alpha_i)$ is constant.

The set of
{\em normalized} \  Seifert invariants $((\alpha_i,\omega_i)_{i=1}^\nu)$
are defined as follows. Write
\begin{equation*}
e=b+\sum \omega_i/\alpha_i
\tag{2}
\end{equation*}
for some integer $b$, and $0\leq \omega_i<\alpha_i$ with
$\omega_i\equiv  -\beta_i \ (\mmod \ \alpha_i)$.  Clearly, these
properties define $\{\omega_i\}_i$ uniquely. Notice that $b\leq
e<0$. For the uniformity of the notations, in the sequel we assume
$\nu\geq 3$. (Recall that for cyclic
quotient singularities Conjecture  was verified in \cite{NN}.)

For each $i$, consider the continued fraction $\alpha_i/ \omega_i=
b_{i1}-1/(b_{i2}-1/(\cdots
-1/b_{i\nu_i})\cdots)$. Then (a
possible) plumbing graph  of $M$ is a star-shaped graph  with
$\nu$ arms. The central vertex has decoration $b$ and the arm
corresponding to
 the index $i$ has $\nu_i$ vertices, and they are
decorated by $b_{i1},\ldots, b_{i\nu_i}$ (the vertex decorated by
$b_{i1}$ is connected by the central  vertex).

We will distinguish those vertices $v\in {\mathcal V}$  of the
graph which have $\delta_v\not= 2$. We will denote by $\bar{v}_0$
the central vertex (with $\delta=\nu$), and by $\bar{v}_i$ the
end-vertex of the $i^{th}$ arm (with $\delta=1$) for all $1\leq
i\leq \nu$.
In this notation, $g_{\bar{v}_0}$ is exactly the
class of the generic fiber. The group $H$ has the following representation:
\begin{equation*}
H=\mbox{ab}\langle g_{\bar{v}_0}, g_{\bar{v}_1},\ldots g_{\bar{v}_\nu}\ |\
g_{\bar{v}_0}^{-b}=\prod_{i=1}^\nu g_{\bar{v}_i}^{\omega_i},
 \ \ g_{\bar{v}_0}=g_{\bar{v}_i}^{\alpha_i}
\  \mbox{for all $i$} \rangle.
\tag{3}
\end{equation*}

 Let $\alpha:=\mbox{lcm}(\alpha_1,\ldots,\alpha_n)$. The order of the
group $H$ and the order $\x $ of the subgroup $\langle g_{\bar{v}_0}
 \rangle$ can be determined by (cf. \cite{Neu}):
\begin{equation*}
|H|=\alpha_1\cdots \alpha_\nu |e|, \ \ \ \ \x:=|\langle
g_{\bar{v}_0} \rangle|=\alpha|e|.
\tag{4}
\end{equation*}

\subsection{Invariants computed from the plumbing graph}\label{ss: inv}
In the sequel we will also  use the  Dedekind  sums. They are
defined as follows \cite{Ra,RG}.
Let $\lfloor x\rfloor$ be  the integer part, and
$\{x\}:= x-\lfloor x\rfloor$ the fractional part of $x$. Then
\[
\bms(h,k)=\sum_{\mu=0}^{k-1}\Big(\Big( \frac{\mu}{k}\Big)\Big)
\Big(\Big( \frac{h\mu}{k} \Big)\Big),\]
where  $((x))$ denotes the Dedekind symbol
\[
((x))=\left\{
\begin{array}{ccc}
\{x\} -1/2 & {\rm if} & x\in {\bR}\setminus {\bZ}\\
0 & {\rm if} & x\in {\bZ}.
\end{array}
\right.
\]
Assume that $M$ is a  Seifert manifold with $e<0$. Then one has the following
formulae for its invariants.

\vspace{2mm}

\noindent $\bullet$ {\bf The Casson-Walker invariant.} \cite{Lescop} (6.1.1):
\begin{equation*}
\frac{24}{|H|}\lambda(M)=
\frac{1}{e}\Big(2-\nu +\sum_{i=1}^\nu\frac{1}{\alpha_i^2}\Big)+
e+3+12\sum_{i=1}^\nu \,  \bms(\beta_i,\alpha_i).
\tag{5}\end{equation*}

\noindent $\bullet$ {\bf $K^2+\#{\mathcal V}$.} \cite{NN} (5.4):
\begin{equation*}
K^2+\#{\mathcal V}=
\frac{1}{e}\Big(2-\nu +\sum_{i=1}^\nu \frac{1}{\alpha_i}\Big)^2 +e +
5+12\sum_{i=1}^\nu \, \bms(\beta_i,\alpha_i).
\tag{6}
\end{equation*}

\noindent $\bullet$ {\bf The Reidemeister-Turaev sign-refined torsion.}
For any $\chi\in \hat{H}$ (and free variable $t\in \C$) set
\begin{equation*}
\hat{P}_{\chi}(t):=
\frac{\bigl(t^{\alpha}\chi(g_{\bar{v}_0})-1\bigr)^{\nu-2}}
{\prod_{i=1}^\nu\ \bigl(t^{\alpha/\alpha_i}\,
\chi(g_{\bar{v}_i})-1\bigr)}.
\tag{7}\end{equation*}
Then, by \cite{NN} (5.8),
the Fourier transform $\hat{\et}_{M,\si_{can}}$ of $\et_{M,\si_{can}}$
is given by
\begin{equation*}
\hat{\et}_{M,\si_{can}}(\bar{\chi})=
\lim_{t\to 1}\, \hat{P}_{\chi}(t)\ \ \
\mbox{for any $\chi\in \hat{H}\setminus \{1\}$}.
\tag{8}
\end{equation*}

\noindent $\bullet$ {\bf The geometric genus of $(X,0)$.} Let
$M$ be a Seifert manifold with $e<0$ and Seifert invariants as above.
Define the Dolgachev-Pinkham (topological) invariant of $M$ by
\begin{equation*}
DP_M:=\sum_{l\geq 0}\, \max\big( 0\, , \, -1+lb-\sum_{i=1}^\nu
\Big\lfloor\, \frac{-l\omega_i}{\alpha_i}\Big\rfloor\, \big).
\tag{9}
\end{equation*}
Assume that $(X,0)$ is a
normal surface singularity  with a good $\C^*$--action (see e.g.
\cite{Pinkh})  such that its link $M$  is a rational homology
sphere.  Then, by \cite{Pinkh}, (5.7):
\begin{equation*}
p_g(X,0)=DP_M.
\tag{10}
\end{equation*}

\section{The main result}

\noindent The key identity of this article is presented in the following
theorem.

\subsection{Theorem.}\label{m1} {\em Let $M$ be  a Seifert 3--manifold with
$e<0$. Then the invariants $\et_{M,\si_{can}}(1)$, $\lambda(M)$,
$K^2+\#{\mathcal V}$ and $DP_M$ are connected by the following identity:}
$$\et_{M,\si_{can}}(1)+\frac{\lambda(M)}{|H|}=\frac{K^2+\#{\mathcal V}}{8}
+DP_M.$$

\subsection{Remark.}\label{mlb} Using \ref{ss: 3.4}(1), for the
(modified) Seiberg-Witten invariant one obtains:
$$\ssw^0_M(\si_{can})=\frac{K^2+\#{\mathcal V}}{8}+DP_M.$$
If $M$ is a singularity link,
one can define on $M$  a canonical {\em contact structure}
$\xi_{can}$ (induced by the natural almost complex structure on
$TM\oplus \R_M$; for  details, see \cite[p. 420]{GS}
or \cite[(4.8)]{NN} ) with $c_1(\xi_{can})$ torsion. On the other hand,
 in \cite{Gompf}, Gompf  associates with such a
 contact structure $\xi$ an invariant
$\theta_M(\xi)$. It turns out that in our case (see \cite[4.8]{NN})
$\theta_M(\xi_{can})=K^2+\#{\mathcal V}-2$. Therefore, for any link of
 singularity which is a Seifert 3-manifold one has:
$$\ssw^0_M(\si_{can})=\frac{\theta_M(\xi_{can})+2}{8}+DP_M.$$

\vspace{2mm}

\noindent The proof of  \ref{m1} is carried out in several steps.

\subsection{Proposition.}\label{m2}
\[
\frac{1}{|H|}\sum_{\chi\in \hat{H}}\, \hat{P}_{\chi}(t)=\sum_{l\geq 0}\, \max
\big( \, 0\, ,\, 1-lb+ \sum_{i=1}^\nu
\Big\lfloor\, \frac{-l\omega_i}{\alpha_i}\Big\rfloor\, \big)\ t^{\x l}.\]
{\em In particular, }
$$\et_{M,\si_{can}}(1)=
\lim_{t\to 1}\,\Big(
\sum_{l\geq 0}\, \max \big(\, 0 \,, \,  1-lb+ \sum_{i=1}^\nu
\Big\lfloor\, \frac{-l\omega_i}{\alpha_i}\Big\rfloor\, \big)\ t^{\x l}
-\frac{1}{|H|}\cdot \hat{P}_1(t)\, \Big).$$
Before we start the proof, we draw  to the reader's attention the
``mysterious'' similarity between our formula (7) for the (Fourier
transform) of the Reidemeister-Turaev torsion, and the formula
\cite[4.2]{Neu} of
W. Neumann of the Poincar\'e series of the graded affine ring associated
with the universal abelian cover of $(X,0)$.

In fact, the next proof is based completely
on Neumann's computation about this
graded ring on [loc. cit.], page 241; (he attributed the idea to D. Zagier).
\begin{proof}
Using the identity $g_{\bar{v}_i}^{\alpha_i}=g_{\bar{v}_0}$ in $H$
(cf. \ref{ss: s1}(3)), first write $\hat{P}_{\chi}(t)$ as
$$(1-t^\alpha\chi(g_{\bar{v}_0}))^{-2}\prod_i\frac{
1-\big(\, t^{\alpha/\alpha_i}\chi(g_{\bar{v}_i})\big)^{\alpha_i}}{
1-t^{\alpha/\alpha_i}\,
\chi(g_{\bar{v}_i})}=
\Big(\, \sum_{s_0=0}^{\infty}\, (1+s_0)\chi(g_{\bar{v}_0})^{s_0}\, t^{\alpha
s_0}\Big)\cdot \prod_i\, \sum_{s_i=0}^{\alpha_i-1} t^{s_i\alpha/\alpha_i}\big(
\chi(g_{\bar{v}_i})\big)^{s_i}
$$
$$=\sum\, (1+s_0)\, t^{\alpha s_0+\sum_i \alpha s_i/\alpha_i}\
\chi\big(g_{\bar{v}_0}^{s_0}g_{\bar{v}_1}^{s_1}\cdots g_{\bar{v}_\nu}^{s_\nu}
\big),$$
where the (unmarked) sum is over $s_0\geq 0$ and
$0\leq s_i< \alpha_i$ for each $i$.
But $\sum_{\chi\in \hat{H}}(h)$ is non-zero
only if $h=1$, and in that case it is $|H|$. Using the group structure
\ref{ss: s1}(3) one gets that all the relations in $H$ have the form
$$
g_{\bar{v}_0}^{l_1+\cdots +l_\nu-lb}\prod_i
g_{\bar{v}_i}^{-\omega_il-\alpha_il_i}=1,$$
where $l_1,\ldots, l_\nu$ and $l$ are integers.
Therefore, $
g_{\bar{v}_0}^{s_0}g_{\bar{v}_1}^{s_1}\cdots g_{\bar{v}_\nu}^{s_\nu}=1$
if and only if
$s_0=l_1+\cdots +l_\nu-lb$ and $s_i=-\omega_il-\alpha_il_i$
\ ($1\leq i\leq \nu$)
for some integers  $l_1,\ldots, l_\nu, l$. Since $0\leq s_i<\alpha_i$
one obtains that
$$l_i=
\Big\lfloor\, \frac{-l\omega_i}{\alpha_i}\, \Big\rfloor.$$
In particular,
$$1+s_0=1-lb +\sum_i
\Big\lfloor\, \frac{-l\omega_i}{\alpha_i}\, \Big\rfloor,$$
and only those integers  $l$ are allowed for which
this number $1+s_0$ is $\geq 1$. It is easy to see that this cannot
happen for $l<0$. Indeed, for $l<0$
$$-lb +\sum_i
\Big\lfloor\, \frac{-l\omega_i}{\alpha_i}\, \Big\rfloor\leq
-lb +\sum_i \frac{-l\omega_i}{\alpha_i}=-le<0.$$
Finally notice that the exponent $\alpha (s_0+\sum_i s_i/\alpha_i)$ of $t$
is $-l\alpha e=l\x$ by \ref{ss: s1}(4).
This ends the proof of the first formula. For the second part, recall that
$\hat{\et}_{M,\si}(1)=0$ (cf. \ref{ss: 3.4}), hence by \ref{ss: inv}(8),
$$\et_{M,\si_{can}}(1)=\lim_{t\to 1}\frac{1}{|H|}
\sum_{\chi\in\hat{H}\setminus \{1\}}\hat{P}_{\chi}(t).$$
\end{proof}
\subsection{Remark.}\label{remark}
Let $(X,0)$ be a normal singularity with a good $\C^*$-action
and affine graded coordinate ring $A=\oplus_kA_k$.
Then its Poincar\'e series is defined by $p_{(X,0)}(t)=\sum_k \dim(A_k)t^k$.
For such an $(X,0)$, the expression from \ref{m2}  provides exactly
$p_{(X,0)}(t^o)$ (cf. \cite{Neu}, page 241). More precisely:
$$p_{(X,0)}(t)=\sum_{l\geq 0}\, \max
\big( \, 0\, ,\, 1-lb+ \sum_{i=1}^\nu
\Big\lfloor\, \frac{-l\omega_i}{\alpha_i}\Big\rfloor\, \big)\ t^{l}.$$
Moreover, if $(X_{ab},0)$ denotes the universal abelian cover of $(X,0)$,
then $p_{(X_{ab},0)}(t)=\hat{P}_1(t)$ (cf. \cite{Neu}, page 240).
Therefore, \ref{m2} reads as follows:
$$\et_{M,\si_{can}}(1)=\lim_{t\to 1}\,\big(
p_{(X,0)}(t^o)-p_{(X_{ab},0)}(t)/|H| \, \big).$$
Notice that for many special  families, the Poincar\'e series
$p_{(X,0)}(t)$  is computed very explicitly, see e.g. \cite{Wag}.

\subsection{Corollary.}\label{m3}
$$\et_{M,\si_{can}}(1)-DP_M=
\lim_{t\to 1}\,\Big(
\sum_{l\geq 0}\, \big( \, 1-lb+ \sum_{i=1}^\nu
\Big\lfloor\, \frac{-l\omega_i}{\alpha_i}\Big\rfloor\, \big)\ t^{\x l}
-\frac{1}{|H|}\cdot \hat{P}_1(t)\, \Big).$$
\begin{proof}
Use \ref{ss: inv}(9), \ref{m2} and the identity $\max(0,x)-\max(0,-x)=x$.
\end{proof}
On the right hand side we have a difference of two series, both having
poles of order two at $t=1$. The next results provide their Laurent
series at $t=1$. In fact, we prefer to expand the series in terms
of the powers of $t^\x-1$ (instead of $t-1$).

\subsection{Proposition.}\label{m4} {\em Define $\chi_M:=2-\sum_{i=1}^\nu\,
(\alpha_i-1)/\alpha_i$ (cf. e.g. with \cite{Neu}). Then
$$\sum_{l\geq 0}\, \big( \, 1-lb+ \sum_{i=1}^\nu
\Big\lfloor\, \frac{-l\omega_i}{\alpha_i}\Big\rfloor\, \big)\ t^{\x l}=
\frac{-e}{(t^\x-1)^2}+\frac{-e-\chi_M/2}{t^\x-1}+
\frac{2-\chi_M}{4}+ \sum_{i=1}^\nu\bms(\beta_i,\alpha_i)+
 R(t),$$
with $\lim_{t\to 1}R(t)=0$. }
\begin{proof}
We follow again \cite{Neu}, page 241. The left hand side of \ref{m4},
via \ref{ss: s1}(2) transforms into
$$\sum_{l\geq 0}\, \Big( -le+\frac{\chi_M}{2}\Big)\,t^{\x l}
+ \sum_{i=1}^\nu \Big(
-\Big\{\, \frac{-l\omega_i}{\alpha_i}\Big\}+\frac{\alpha_i-1}{2
\alpha_i}\,  \Big)\, t^{\x l}.$$
Evidently
$$\sum_{l\geq 0}\, \Big( \, -le+\frac{\chi_M}{2}\Big) \ t^{\x l}=
\frac{-et^{\x}}{(1-t^\x)^2}+\frac{\chi_M/2}{1-t^\x},$$
which gives the non-holomorphic part.
The second contribution is a sum over $1\leq i\leq\nu$. For each fixed $i$,
write $l=\alpha_i m+q$ with $m\geq 0$ and $0\leq q<\alpha_i$. Using the
 notation $\sum_q:=\sum_{q=0}^{\alpha_i-1}$ and $\sum_q':=
\sum_{q=1}^{\alpha_i-1}$, the $i^{th}$ summand is
$$\sum_q\Big(
-\Big\{\, \frac{-q\omega_i}{\alpha_i}\Big\}+\frac{\alpha_i-1}{2\alpha_i}\,
  \Big)\ \sum_{m\geq 0}t^{\x \alpha_i m+\x q}=
\frac{\sum_q\big(
-\big\{\, \frac{-q\omega_i}{\alpha_i}\big\}+\frac{\alpha_i-1}{2\alpha_i}\,
  \big)\ t^{\x q}}{1-t^{\x \alpha_i}}.$$
Separating the two cases  $q=0$ and $q>0$, and using
the definition of the Dedekind symbol and the identity
$\{-x\}=1-\{x\}$ for $x\not\in \Z$, this is transformed into
$$A(t):=\frac{
\frac{\alpha_i-1}{2\alpha_i}+\sum_q'\big(\big( \frac{q\omega_i}
{\alpha_i}\big)\big) t^{\x q}-\frac{1}{2\alpha_i}(t^{\x}+t^{2\x}+\cdots
+t^{(\alpha_i-1)\x})}
{1-t^{\x \alpha_i}}.$$
By L'Hospital theorem  (and by  some simplifications):
$$\lim_{t\to 1}A(t)=-{\sum_q}'\Big(\Big(\frac{q\omega_i}{\alpha_i}\Big)\Big)
\frac{q}{\alpha_i}+\frac{\alpha_i-1}{4\alpha_i}.$$
Since $\sum_q'((q\omega_i/\alpha_i))=0$ and $\omega_i\equiv -\beta_i$
 ($\mmod\ \alpha_i$), the result follows from the definition of the
Dedekind symbol and the Dedekind sums.
\end{proof}
\subsection{Remarks.}\label{m5} Cf. \cite{Neu},  page 242. In fact,
$$\sum_{l\geq 0}\, \big( \, 1-lb+ \sum_{i=1}^\nu
\Big\lfloor\, \frac{\omega_i}{\alpha_i}(-l)\Big\rfloor\, \big)\ t^{\x l}=
\frac{-e}{(t^\x-1)^2}+\frac{-e-\chi_M/2}{t^\x-1}+
 \sum_{i=1}^\nu \frac{1}{\alpha_i} {\sum_{\Z_{\alpha_i}}}'
\frac{1}{(1-\xi)(1-\xi^{\omega_i}t^\x)},$$
where the last   sum is over $\xi^{\alpha_i}=1$, $\xi\not=1$.
Since we do not need this statement now,   we will skip its proof.
The interested reader can prove it easily using the expression $A(t)$ above
and the property (16c) of the Dedekind symbol from \cite{RG}, page 14.

\vspace{2mm}

Paradoxically, the Laurent expansion of the elementary rational fraction
 $\hat{P}_1(t)$ is more complicated.

\subsection{Proposition.}\label{m6} {\em
$$\frac{\hat{P}_1(t)}{|H|}=
\frac{-e}{(t^\x-1)^2}+\frac{-e-\chi_M/2}{t^\x-1}+ E+  Q(t),$$
where $\lim_{t\to 1}Q(t)=0$ and }
$$E:=-\frac{(e+1)(e+5)}{12e}+\frac{1}{4}\sum_i\big(1-\frac{1}{\alpha_i}\big)
+\frac{1}{12e}\sum_i\big(1-\frac{1}{\alpha_i}\big)
\big(4+\frac{1}{\alpha_i}\big)-
\frac{1}{4e}\sum_{i<j}\big(1-\frac{1}{\alpha_i}\big)
\big(1-\frac{1}{\alpha_j}\big).
$$
\begin{proof}
First notice that one has the following Taylor expansion;
$$\frac{t^\gamma-1}{t^\tau-1}=\frac{\gamma}{\tau}+\frac{\gamma}
{2\x \tau}( \gamma-\tau)\cdot (t^\x-1)+
\frac{\gamma}{\x \tau}( \gamma-\tau)\Big(
\frac{2\gamma-\tau}{12\x}-\frac{1}{4}\Big) \cdot (t^\x-1)^2+\cdots .$$
Now,  use this formula $\nu+2$ times in the expression
$$\hat{P}_1(t)=\frac{1}{(t^\x-1)^2}\cdot \Big(
\frac{t^\x-1}{t^\alpha-1}\Big)^2\cdot
\prod_i\frac{t^\alpha-1}{t^{\alpha/\alpha_i}-1}.$$
A long (but elementary) computation, together with \ref{ss: s1}(4),
gives the result.
\end{proof}

\vspace{2mm}

\noindent {\em Proof of \ref{m1}.} Apply \ref{ss: inv}(5) and (6),
respectively \ref{m3}, \ref{m4} and \ref{m6}. The verification is
elementary. \qed

\subsection{Remark. The Reidemeister-Turaev torsion revisited.}\label{rt}
In the above results we needed to determine only $\et_{M,\si_{can}}(1)$,
i.e. only that sign refined torsion which is associated with the
canonical $spin^c$ structure.  But, via similar computations, we can obtain
the complete set of invariants $\{\et_{M,\si}(1)\}_{\si\in Spin^c(M)}$.

Indeed, we recall that $Spin^c(M)$ is an $H$-torsor. Let $\si$ be an
arbitrary element of $Spin^c(M)$, and take (the unique) $h_{\si}\in H$ so that
$h_{\si}\cdot \si_{can}=\si$. Then, by \cite[5.7]{NN}, one has:
$$\hat{\et}_{M,\si}(\bar{\chi})=\bar{\chi}(h_{\si})\cdot \lim_{t\to 1}
\hat{P}_{\chi}(t)\ \ \mbox{for any $\chi\in \hat{H}\setminus \{1\}$}.$$
Using the representation \ref{ss: s1}(3) of $H$, $h_\si$ evidently
can be written as
$$h_\si=g_{\bar{v}_0}^{a_0}g_{\bar{v}_1}^{a_1}\cdots g_{\bar{v}_\nu}^{a_\nu}
\ \ \ \mbox{for some integers} \ a_0,a_1,\ldots,a_\nu.$$
Then repeating the arguments from the proof of \ref{m2}, one gets:
$$\bar{\chi}(h_\si)\hat{P}_{\chi}(t)
=\sum\, (1+s_0)\, t^{\alpha s_0+\sum_i \alpha s_i/\alpha_i}\
\chi\big(h_{\si}^{-1}
g_{\bar{v}_0}^{s_0}g_{\bar{v}_1}^{s_1}\cdots g_{\bar{v}_\nu}^{s_\nu}\big),$$
with $s_0\geq 0$ and  $0\leq s_i< \alpha_i$ for each $i$.
But  $h_\si
g_{\bar{v}_0}^{s_0}g_{\bar{v}_1}^{s_1}\cdots g_{\bar{v}_\nu}^{s_\nu}=1$
if and only if
$s_0=a_0+l_1+\cdots +l_\nu-lb$ and $s_i=a_i-\omega_il-\alpha_il_i$
\ ($1\leq i\leq \nu$)
for some integers  $l_1,\ldots, l_\nu, l$. Since $0\leq s_i<\alpha_i$
one obtains that
$$l_i=
\Big\lfloor\, \frac{-l\omega_i+a_i}{\alpha_i}\, \Big\rfloor.$$
In particular,
$$1+s_0=1+a_0-lb +\sum_i
\Big\lfloor\, \frac{-l\omega_i+a_i}{\alpha_i}\, \Big\rfloor,$$
and only those integers  $l$ are allowed for which
 $1+s_0\geq 1$.
Therefore, one gets:

\subsection{Proposition.}\label{tors} {\em Assume that $M$ is a Seifert
3-manifold with $e<0$. Fix an arbitrary $spin^c$ structure
$\si\in Spin^c(M)$  characterized by $h_{\si}\cdot \si_{can}
=\si$. Write $h_\si$ as
$g_{\bar{v}_0}^{a_0}g_{\bar{v}_1}^{a_1}\cdots g_{\bar{v}_\nu}^{a_\nu}$
for some integers $a_0,a_1,\ldots, a_\nu$.
Finally, define
$\tilde{a}:=\alpha\cdot ( a_0+\sum_ia_i/\alpha_i)$. Then
$$\et_{M,\si}(1)=
\lim_{t\to 1}\,\Big(
\sum_{l \in \Z}\, \max \big(\, 0 \,, \,  1+a_0-lb+ \sum_{i=1}^\nu
\Big\lfloor\, \frac{-l\omega_i+a_i}{\alpha_i}\Big\rfloor\, \big)\ t^{\x l+
\tilde{a}}
-\frac{1}{|H|}\cdot \hat{P}_1(t)\, \Big).$$
(Above, there exists an integer $m$, such that for $l<m$ the contribution in
the above  sum is zero.)
}

\end{document}